# CONFORMAL CHARACTERS AND THETA SERIES

Wolfgang Eholzer and Nils-Peter Skoruppa

ABSTRACT. We describe the construction of vector valued modular forms transforming under a given congruence representation of the modular group $\mathrm{SL}(2,\mathbb{Z})$ in terms of theta series. We apply this general setup to obtain closed and easily computable formulas for conformal characters of rational models of $\mathcal{W}$-algebras.

## Contents

## 1. Introduction

Once a rational model of a $\mathcal{W}$-algebra is constructed or its existence is conjectured the next natural step is to compute explicitly its conformal characters. This is interesting for giving further evidence to the existence of the rational model in question, if its existence is not sure, and it is interesting in that conformal characters are usually surprisingly distinguished modular functions which are often related to a priori purely number theoretical objects like generalized Rogers-Ramanujan identities or dilogarithm identities.

However, in general it is quite difficult to construct conformal characters directly from physical information about the $\mathcal{W}$-algebra; even the direct computation of their first few Fourier coefficients affords considerable computer power.

In [E-S] we took a different point of view in order to determine conformal characters amongst all possible modular functions.

We formulated a list of certain axioms which hold true for all known sets of conformal characters of rational models of $\mathcal{W}$-algebras. The only data from an underlying rational model which occur in these axioms are respectively its central charge and its conformal dimensions. We showed that, for several rational models, these axioms uniquely determine the conformal characters belonging to a given list of central charge and conformal dimensions.

Research at MSRI is supported in part by NSF grant no. DMS–9022140

Thus, once the central charge and conformal dimensions of a rational model are known, the computation of its conformal characters can be viewed as a problem which is completely independent from the theory of $\mathcal{W}$-algebras, i.e., for this computation one is left with a construction problem, namely, the problem of finding, by whatever means, a set of functions fulfilling the indicated list of axioms.

The purpose of the present note is to describe such a mean which can solve in many cases this construction problem. In particular, we shall apply our method to the case of five special rational models. The reason for the choice of these models is that the representation theory of the $\mathrm{SL}(2,\mathbb{Z})$-representation on their conformal characters can be treated homogeneously in some generality, and that the conformal characters of one of these models (of type $\mathcal{W}(2,8)$ with central charge $c = -\frac{3164}{23}$) could not be computed by the so far known methods.

This note is organized as follows: In section 2 we give a short description of the five special models. We recall the axiomatic description of conformal characters, and, in particular, we recall the results obtained in [E-S] which concern the five special models. In Section 3 we describe a general procedure for the construction of vector valued modular forms transforming under a given matrix representation of $\mathrm{SL}(2,\mathbb{Z})$. As already mentioned, this procedure is thought to be useful in general for finding explicit and easily computable formulas for conformal characters. In section 4 we apply the general setup of the preceding section to the case of the five special models, and we derive explicit formulas for their conformal characters. Finally, in §5 we compare our results with those formulas for the conformal characters of the five models which can be obtained (assuming certain conjectures) from the representation theory of Casimir-$\mathcal{W}$-algebras.

**Notation.** We use $\mathfrak{H}$ for the complex upper half plane, $\tau$ as a variable in $\mathfrak{H}$, $q = e^{2\pi i \tau}$, $q^\delta = e^{2\pi i \delta \tau}$, $T = \begin{pmatrix} 1 & 1 \\ 0 & 1 \end{pmatrix}$, $S = \begin{pmatrix} 0 & -1 \\ 1 & 0 \end{pmatrix}$, $\Gamma$ for the group $\mathrm{SL}(2,\mathbb{Z})$, and

$$\Gamma(n) = \{A \in \mathrm{SL}(2,\mathbb{Z}) \mid A \equiv \mathrm{id} \pmod{n}\}$$

for the principal congruence subgroup of $\mathrm{SL}(2,\mathbb{Z})$ of level $n$. We use $\eta$ for the Dedekind eta function

$$\eta(\tau) = e^{\pi i \tau / 12} \prod_{n \geq 1} (1 - q^n).$$

The group $\Gamma$ acts on $\mathfrak{H}$ by

$$A\tau = \frac{a\tau + b}{c\tau + d} \qquad \left(A = \begin{pmatrix} a & b \\ c & d \end{pmatrix}\right).$$

For a complex vector valued function $F(\tau)$ on $\mathfrak{H}$, and for an integer $k$ we use $F|_k A$ for the function defined by

$$(F|_k A)(\tau) = (c\tau + d)^{-k} F(A\tau).$$

Finally, for a matrix representation $\rho: \Gamma \to \mathrm{GL}(n,\mathbb{C})$ and an integer $k$ we use $M_k(\rho)$ for the vector space of all holomorphic maps $F: \mathfrak{H} \to \mathbb{C}^n$ (= column vectors) which satisfy $F|_k A = \rho(A) F$ for all $A \in \Gamma$, and which are bounded in any region $\mathrm{Im}(\tau) \geq r > 0$. Thus, if $\rho$ is the trivial representation, then $M_k(\rho)$ is the space of ordinary modular forms on $\Gamma$ and of weight $k$.

## 2. An axiomatic characterization of the conformal characters of five special rational models

In this section we recall some results on the conformal characters of 5 rational models, which were, among others, considered in [E-S]. The types, the central charges, the sets of conformal dimensions $H_c$, and the effective central charge $\tilde{c} = c - 24 \min H_c$ of these models are listed in Table 1; for a more detailed description the reader is referred to loc. cit.. Denote by $\chi_{c,h}$ ($h \in H_c$) the conformal characters of the rational model with central charge $c$, where we choose the notation such that

$$\chi_{c,h} = q^{h-c/24} \cdot (\text{power series in } q).$$

Table 1: Central charges and conformal dimensions

| type | $c$ | $\tilde{c}$ | $H_c$ |
|---|---|---|---|
| $\mathcal{W}_{G_2}(2, 1^{14})$ | $-\frac{8}{5}$ | $\frac{16}{5}$ | $\frac{1}{5}\{0, -1, 1, 2\}$ |
| $\mathcal{W}_{F_4}(2, 1^{26})$ | $\frac{4}{5}$ | $\frac{28}{5}$ | $\frac{1}{5}\{0, -1, 2, 3\}$ |
| $\mathcal{W}(2, 4)$ | $-\frac{444}{11}$ | $\frac{12}{11}$ | $-\frac{1}{11}\{0, 9, 10, 12, 14, 15, 16, 17, 18, 19\}$ |
| $\mathcal{W}(2, 6)$ | $-\frac{1420}{17}$ | $\frac{20}{17}$ | $-\frac{1}{17}\{0, 27, 30, 37, 39, 46, 48, 49, 50,$ $52, 53, 55, 57, 58, 59, 60\}$ |
| $\mathcal{W}(2, 8)$ | $-\frac{3164}{23}$ | $\frac{28}{23}$ | $-\frac{1}{23}\{0, 54, 67, 81, 91, 94, 98, 103, 111,$ $112, 116, 118, 119, 120, 122, 124,$ $125, 129, 130, 131, 132, 133\}$ |

Fix now one of the central charges $c$ of Table 1. Suppose that we are given functions $\xi_{c,h}$ ($h \in H_c$), defined on the upper half plane $\mathfrak{H}$, which satisfy the following list of properties:

**Properties of conformal characters.**
(1) The functions $\xi_{c,h}$ are nonzero modular functions for some congruence subgroup of $\Gamma = \mathrm{SL}(2, \mathbb{Z})$ without poles in the upper half plane.
(2) The space of functions spanned by the $\xi_{c,h}$ ($h \in H_c$) is invariant under $\Gamma$ with respect to the action $(A, \xi) \mapsto \xi(A\tau)$.
(3) For each $h \in H_c$ one has $\xi_{c,h} = \mathcal{O}(q^{-\tilde{c}/24})$ as $\mathrm{Im}(\tau)$ tends to infinity, where $\tilde{c} = c - 24 \min H_c$.
(4) For each $h \in H_c$ the function $q^{-(h-\frac{c}{24})}\xi_{c,h}$ is periodic with period 1.
(5) The Fourier coefficients of the $\xi_{c,h}$ are rational numbers.

It was shown in [E-S] that the functions $\xi_{c,h}$ are uniquely determined by these five properties, up to multiplication by scalars. Very likely the conformal characters $\chi_{c,h}$ satisfy the properties (1) to (5) if we view them as functions of $\tau$ by setting $q = \exp(2\pi i \tau)$ (for a detailed discussion of this conjecture cf. loc. cit.). Thus, assuming this conjecture for a moment, we conclude that, for any set of functions $\xi_{c,h}$ as above, we have $\chi_{c,h} = \mathrm{const.} \cdot \xi_{c,h}$.

As already explained the main purpose of this paper is the explicit construction of functions $\xi_{c,h}$ satisfying (1) to (5). To this end we have to recall more precisely what we showed in [E-S].

Denote by $\xi_c$ the column vector whose entries are the $\xi_{c,h}$ ordered in some way. Let $l$ be the denominator of $c$, let $k$ be the integer associated to $c$ by Table 2, and let $\rho_l$ be the $\Gamma$-representation defined by

$$(\eta^{2k}\xi)|_k A = \rho_l(A)\, \eta^{2k}\xi \qquad (A \in \Gamma).$$

As we shall see in a moment the a priori two different representations with $l = 5$ are in fact equivalent; so for notational convenience we denote them by the same symbols. We showed in [E-S; cf. §4.4] that properties (1) to (5) imply

**Theorem (Characterization of $\rho_l$).** *The representation $\rho_l$ is irreducible, its kernel contains $\Gamma(l)$, and it takes its values in $\mathrm{GL}(l-1, \mathbb{Q}(e^{2\pi i/l}))$.*

Consulting any table of irreducible representations of $\Gamma/\Gamma(l) \approx \mathrm{SL}(2, \mathbb{Z}/l\mathbb{Z})$ one verifies that the three properties listed in the Theorem uniquely determine $\rho_l$ up to equivalence. We shall give an explicit description of $\rho_l$ in § 4.

From property (3) we have $\eta^{2k}\xi_c = \mathcal{O}(q^\delta)$ for $q \to 0$, where $\delta = -\tilde{c} + k/12$, and, in particular, that $\eta^{2k}\xi_c$ is an element of $M_k(\rho_l)$. The dimensions of these spaces can be computed using the dimension formula in [E-S]. The resulting dimensions and the values of $\delta$ are listed in Table 2.

Let $M_k^{(\delta)}(\rho_l)$ be the subspace of all $F \in M_k(\rho_l)$ satisfying $F = \mathcal{O}(q^\delta)$. In [E-S] it was shown that this subspace is one-dimensional, which, by obvious arguments, implies that $\xi_c$ is unique up to multiplication by diagonal matrices (actually, it was shown in loc. cit. that $M_h(\rho_l \otimes \theta^{2h-2k})$ is one-dimensional, where $\theta^2(A) = (\eta^2|_1 A)/\eta^2$. However, this latter space is obviously isomorphic to $M_k^{(\delta)}(\rho_l)$ via multiplication by $\eta^{2k-2h}$).

Table 2: Certain data related to the five rational models

| $\mathcal{W}$-algebra | $c$ | $l$ | $k$ | $\delta$ | $\dim M_k(\rho_l)$ |
|---|---|---|---|---|---|
| $\mathcal{W}_{G_2}(2, 1^{14})$ | $-\frac{8}{5}$ | 5 | 4 | $\frac{1}{5}$ | 1 |
| $\mathcal{W}_{F_4}(2, 1^{26})$ | $\frac{4}{5}$ | 5 | 10 | $\frac{3}{5}$ | 3 |
| $\mathcal{W}(2, 4)$ | $-\frac{444}{11}$ | 11 | 6 | $\frac{5}{11}$ | 5 |
| $\mathcal{W}(2, 6)$ | $-\frac{1420}{17}$ | 17 | 2 | $\frac{2}{17}$ | 2 |
| $\mathcal{W}(2, 8)$ | $-\frac{3164}{23}$ | 23 | 10 | $\frac{18}{23}$ | 17 |

3. Realization of modular representations by theta-series

In this section we show how one can, under certain hypothesis, construct systematically vector valued modular in $M_k(\rho)$ for a given matrix representation $\rho$ of $\Gamma$ and given weight $k$.

The first step is a realization of $\rho$ as subrepresentation of a Weil representation. We explain this notion: Let $M$ be a finite abelian group and

$$\mathcal{Q}\colon M \to \mathbb{Q}/\mathbb{Z}$$

a non-degenerate quadratic form. By quadratic form we mean that $\mathcal{Q}$ is even (i.e. $\mathcal{Q}(-x) = \mathcal{Q}(x)$) and the map

$$(x, y) \mapsto \mathcal{B}(x, y) := \mathcal{Q}(x + y) - \mathcal{Q}(x) - \mathcal{Q}(y)$$

is $\mathbb{Z}$-bilinear, and by non-degenerate we mean that for each $x \not= 0$ there exists a $y$ such that $\mathcal{B}(x,y) \not= 0$. We shall call such a pair $(M, \mathcal{Q})$ a quadratic module. The Weil representation $\omega = \omega_{(M,\mathcal{Q})}$ associated to $M$ and $\mathcal{Q}$ is a projective right-representation of $\Gamma$ on the space $\mathbb{C}^M$ of complex-valued functions $f$ defined on $M$. The operators corresponding under $\omega$ to the generators $S$ and $T$ of $\Gamma$ are given by

$$f|\omega(T)(x) = \mathrm{e}(\mathcal{Q}(x))\,f(x), \quad f|\omega(S)(x) = \gamma(1)^{-1} \sum_{y \in M} \mathrm{e}(-\mathcal{B}(x,y))\,f(y).$$

Here, for any integer $a$ we use

$$\gamma(a) = \sum_{x \in M} \mathrm{e}(a\mathcal{Q}(x)).$$

Moreover, here and in the following, for any $x \in \mathbb{Q}/\mathbb{Z}$, we shall write $\mathrm{e}(x)$ for $\mathrm{e}^{2\pi i r}$, where $r$ is any representative of $x$. The reader is referred to [N] for more details on Weil representations (To translate our terminology to the one in [N] note that $f|\omega_{(M,\mathcal{Q})}(A) = w_{(M,-\mathcal{Q})}(A^{-1})f$, where $w_{(M,-\mathcal{Q})}$ is the Weil (left-)representation associated to the quadratic module $(M, -\mathcal{Q})$ as in [N].)

Denote by $l$ the level of $\mathcal{Q}$, i.e. the smallest positive integer such that $l\mathcal{Q}(M) = 0$. Using that $\mathcal{Q}$ is non-degenerate it is easy to show that $|\gamma(a)| = |M|^{1/2}$ for any integer relatively prime to $l$. Moreover, $\gamma(a)$ is obviously an element of $K = \mathbb{Q}(\mathrm{e}^{2\pi i/l})$. We identify $(\mathbb{Z}/l\mathbb{Z})^*$ with the Galois group $G$ of $K$ in the usual way such that $a \in (\mathbb{Z}/l\mathbb{Z})^*$ corresponds to the Galois substitution $\sigma_a$ which is given by $\mathrm{e}(1/l)^{\sigma_a} = \mathrm{e}(a/l)$. For $a \in G$ set $co(a) = \gamma(a)/\gamma(1)$, and denote by $\mu_8$ the group of eights roots of unity. Then it is easy to show that $a \mapsto co(a)$ defines a co-cycle of $G$ with values in $\mu_8$, i.e. $co(a)$ is a eights root of unity and one has

$$co(a)^{\sigma_b} co(b) = co(ab) \qquad (a, b \in (\mathbb{Z}/l\mathbb{Z})^*).$$

To verify this note that $\gamma(a) = \gamma(1)^{\sigma_a}$, from which the co-cycle relation is obvious. Moreover, this identity (together with $|\gamma(a)| = |M|^{1/2}$) implies $|co(a)^\sigma| = 1$ for all $\sigma \in G$, i.e. that $co(a)$ is a root of unity, say $co(a) = \mathrm{e}(r/l)$. Since $co(a)$ is obviously invariant under $G^2$ we find that $r(b^2 - 1)$, for any $b$ relatively prime to $l$, is divisible by $l$, which is only possible if $l/(l,r)$ divides 8, i.e. if $co(a) \in \mu_8$.

The Weil representation $\omega$ is a true representation, i.e. not only a projective one, if and only if the co-cycle $co$ is a character [N, Satz 1]. If this is the case we call $\omega$ a proper Weil representation. From the discussion above it is clear that $co$ is a character if and only if $co$ takes values in $\{\pm 1\}$. Moreover, if $\omega$ is proper, then its kernel contains $\Gamma(l)$ [N, Satz 2].

One now has the following theorem.

**Theorem** [N-W]. *Each irreducible right-representation of $\Gamma$ whose kernel contains a principal congruence subgroup is isomorphic to a subrepresentation of a suitable proper Weil representation.*

We call two quadratic modules $(M, \mathcal{Q})$ and $(M', \mathcal{Q}')$ isomorphic if there exist an isomorphism (of abelian groups) $\pi \colon M \to M'$ such that $\mathcal{Q}' \circ \pi = \mathcal{Q}$, and we denote such an isomorphism by

$$\pi \colon (M, \mathcal{Q}) \xrightarrow{\approx} (M', \mathcal{Q}').$$

It is easy to show that isomorphic quadratic modules yield isomorphic Weil representations: an isomorphism of (projective or proper) $\Gamma$-representations is given by the map
$$\pi^*: \mathbb{C}^{M'} \to \mathbb{C}^M, \qquad f \mapsto \pi^* f = f \circ \pi.$$

As the next step for constructing elements of spaces $M_k(\rho)$ we connect Weil representations and theta series by lifting quadratic modules to lattices and quadratic forms on them.

More precisely, let $(M, \mathcal{Q})$ be a quadratic module. Assume that $L$ is a complete lattice in some rational finite-dimensional vector space $V$ and $Q$ a positive definite non-degenerate quadratic form on $V$ which takes on integral values on $L$, and such that there exists an isomorphism of quadratic modules
$$\pi: (L^\sharp/L, \widetilde{Q}) \xrightarrow{\approx} (M, \mathcal{Q}).$$

Here we use $L^\sharp$ for the dual lattice of $L$ with respect to $Q$, i.e. $L^\sharp$ is the set of all $y \in V$ such that $B(L, y) \subset \mathbb{Z}\}$ with $B(x,y) = Q(x+y) - Q(x) - Q(y)$, and we use $\widetilde{Q}$ for the induced quadratic form
$$\widetilde{Q}: L^\sharp/L \to \mathbb{Q}/\mathbb{Z}, \qquad x + L \mapsto Q(x) + \mathbb{Z}.$$

We shall call such a pair $(L, Q)$ a lift of the quadratic module $(M, \mathcal{Q})$.

Let $p$ a homogeneous spherical polynomial on $V$ with respect to $Q$ of degree $\nu$, i.e. if we choose a basis $b_j$ of $V$, then $p\left(\sum b_j \xi_j\right)$ becomes a complex homogeneous polynomial in the variables $\xi_j$ of degree $\nu$ satisfying
$$\nabla G^{-1} \nabla' p\Big(\sum_j b_j \xi_j\Big) = 0,$$
where $\nabla = (\frac{\partial}{\partial \xi_1}, \dots)$ and $G = (B(b_j, b_k))_{j,k}$ is the Gram matrix of $B$.

Finally, for $f \in \mathbb{C}^M$, set
$$\theta_f = \sum_{x \in L^\sharp} (\pi^* f)(x)\, p(x)\, q^{Q(x)}.$$

Here we view $\pi^* f$ as function on function on $L^\sharp$ which is periodic with period lattice $L$.

We assume that $V$ has even dimension $2r$. Then the Weil representation $\omega = \omega_{(M,\mathcal{Q})}$ is proper as follows from the discussion of the co-cycle $co(a)$ above and Milgram's theorem which implies that $co(a)$ is $\pm 1$ (cf. Appendix below). One has

**Theorem (Representation by theta series).** *The map $\mathbb{C}^M \ni f \mapsto \theta_f$ has the property $\theta_f|_{r+\nu} A = \theta_{f|\omega(A)}$ for all $A \in \Gamma$, i.e. it defines a homomorphism of $\Gamma$-modules.*

This is, in various different formulations, a well-known theorem. For the reader's convenience we shall sketch the proof in the appendix below.

Let now $\rho: \Gamma \to \mathrm{GL}(n, \mathbb{C})$ be a congruence matrix representation, and assume that we have determined a quadratic module $(M, \mathcal{Q})$ such that the associated Weil representation is proper and contains a subrepresentation which is isomorphic to the (right-)representation $\mathbb{C}^n \times \Gamma \ni (z, A) \mapsto z\rho(A)'$, where the prime denotes transposition. The existence of such a $(M, \mathcal{Q})$ is guaranteed by the first theorem. Thus, there exists a $\Gamma$-invariant subspace of $\mathbb{C}^M$ with basis $f_j$ such that

$$\Phi|\omega(A) = \rho(A)\Phi \qquad (A \in \Gamma),$$

where $\Phi$ denotes the column vector build from the $f_j$. Assume furthermore that there exists a lift $(L, Q)$ of $(M, \mathcal{Q})$, i.e. an isomorphism

$$\pi\colon (L^\sharp/L, \widetilde{Q}) \xrightarrow{\approx} (M, \mathcal{Q})$$

with a lattice $L$ of even rank $2r$. Let $p$ be a homogeneous spherical polynomial w.r.t. $Q$ of degree $\nu$. From the last theorem it is then clear that we have the following

**Theorem (Realization by theta series).** *The function*

$$\theta = \sum_{x \in L^\sharp} \Phi(\pi(x))\, p(x)\, q^{Q(x)}$$

*is an element of $M_{r+\nu}(\rho)$.*

**Appendix.** We proof the theorem on representation by theta series. It is consequence of the following

**Lemma (Basic Transformation formula).** *Let $L$ be a lattice in a rational vector space $V$ of dimension $2r$, let $Q$ be a positive definite quadratic form on $V$ which takes on integral values on $L$, let $L^\sharp$ and $B$ be defined as above, let $w \in V \otimes \mathbb{C}$ with $Q(w) = 0$, let $\nu$ a non-negative integer, and let $z \in V$. Then one has*

$$\tau^{-r-\nu} \sum_{x \in L} [B(w, x+z)]^\nu\, \mathrm{e}(-Q(x+z)/\tau)$$

$$= \frac{i^{-r}}{\sqrt{[L^\sharp : L]}} \sum_{y \in L^\sharp} [B(w, y)]^\nu\, \mathrm{e}(\tau Q(y)^t - B(y, z)),$$

*where $\tau$ is a variable in the complex upper half plane.*

The lemma is a well-known consequence of the Poisson summation formula; for a proof cf. [Sch, p. 206]. (For verifying that our formula is equivalent to the one given loc. cit. identify $L$ with $\mathbb{Z}^{2r}$ by choosing a $\mathbb{Z}$-basis $b_j$ of $L$, and note that then $L^\sharp = G^{-1}\mathbb{Z}^{2r}$ and $\det(G) = [L^\sharp : L]$ where $G = (B(b_j, b_k))$ is the Gram matrix of $L$. Moreover, the transformation formula loc. cit. is only stated for $\tau = it$ ($t$ real); the general formula follows by analytic continuation.)

*Proof of the theorem on representation by theta series.* Since any homogeneous spherical polynomial of degree $\nu$ can be written as linear combination of the special ones $B(x, w)^\nu$ (where $w \in \mathbb{C}$, $Q(w) = 0$) we can assume that $p$ is of this special form. Since $S$ and $T$ generate $\Gamma$ it suffices to prove the asserted formula for these elements. For $A = T$ the formula is obvious. For proving the case $A = S$ let in the basic transformation formula $z$ be an element of $L^\sharp$, multiply by $f(z)$ and sum over a set of representatives $z$ for $L^\sharp/L$. Using

$$\sum_{x \in L^\sharp/L} \mathrm{e}(Q(x)) = i^r \sqrt{[L^\sharp : L]}$$

(Milgram's theorem, e.g. [H-M, p. 127]) we realize the claimed formula. □

## 4. Theta series associated to quaternion algebras and the conformal characters of the five special model

We now follow the procedure outlined in the foregoing section to construct elements of $M_k(\rho_l)$ where $l$ denotes an odd prime $l \equiv -1 \bmod 3$ and $\rho_l$ is the matrix representation introduced in §2.

We first describe how to obtain $\rho_l$ from a proper Weil representation. Let $\omega$ be the Weil representation associated to the quadratic module $(\mathbb{F}(l^2), \mathrm{n}(x)/l)$. Here $\mathbb{F}(l^2)$ is the field with $l^2$ elements, and $\mathrm{n}(x) = x \cdot \overline{x}$ with $x \mapsto \overline{x} = x^l$ denoting the non-trivial automorphism of $\mathbb{F}(l^2)$. Note that $\mathrm{tr}(x\overline{y})/l$ where $\mathrm{tr}(x) = x + \overline{x}$ is the bilinear form associated to $n(x)/l$. The Weil representation $\omega$ associated is thus a (right-)representation of $\Gamma$ on the space of functions $f \colon \mathbb{F}(l^2) \to \mathbb{C}$, and it is given by

$$f|\omega(T)(x) = \mathrm{e}(\mathrm{n}(x)/l)\, f(x), \qquad f|\omega(S)(x) = \frac{-1}{l} \sum_{y \in \mathbb{F}(l^2)} \mathrm{e}(-\mathrm{tr}(\overline{x}y)/l)\, f(y).$$

Here we used

$$\sum_{x \in \mathbb{F}(l^2)} \mathrm{e}(\mathrm{n}(x)/l) = -l,$$

as follows for instance from Milgram's theorem and the considerations below where we shall obtain $\mathbb{F}(l^2) = L^\sharp/L$ with a lattice of rank 4. Note that this identity implies in particular that $\omega$ is a proper representation (cf. the discussion in the preceding section).

Let $\chi$ be one of the two characters of order 3 of the multiplicative group of nonzero elements in $\mathbb{F}(l^2)$, and let $G$ be the subgroup of elements with $\mathrm{n}(x) = 1$. Note that the existence of $\chi$ follows from the assumption $l \equiv -1 \bmod 3$. Let $X(\chi)$ be the subspace of all $\phi \in X$ which satisfy $\phi(gx) = \chi(g)\phi(x)$ for all $g \in G$. It is easily checked that $X(\chi)$ is a $\Gamma$-submodule of $X$. In fact, it is even an irreducible one [N-W, Satz 2]. As basis for $X(\chi)$ we may pick the functions $\chi_r$ $(1 \le r \le l-1)$ which are defined by $\chi_r(x) = \chi(x)$ if $\mathrm{n}(x) = r$ and $\chi_r(x) = 0$ otherwise. Let $\Phi_\chi$ be the complex column vector valued function on $\mathbb{F}(l^2)$ whose $r$-th component equals $\chi_r$. We then have $\Phi_\chi|A = \rho(A)\Phi_\chi$ with a unique matrix representation $\rho \colon \Gamma \to \mathrm{GL}(l-1, \mathbb{C})$. It is an easy exercise to verify the identities

$$\rho(T) = \mathrm{diag}(\mathrm{e}^{2\pi i 1/l}, \cdots, \mathrm{e}^{2\pi i (l-1)/l}), \quad \rho(S) = (\lambda(rs))_{1 \le r, s \le l-1},$$

where we use

$$\lambda(r) = \frac{-1}{l} \sum_{\substack{x \in \mathbb{F}(l^2) \\ \mathrm{n}(x) = r}} \chi(x)\, \mathrm{e}(\mathrm{tr}(x)/l).$$

(In the identity $\mathrm{n}(x) = r$ the $r$ has to be viewed as an element of $\mathbb{F}(l^2)$.)

Note that $\lambda(r)$ does not depend on the choice of $\chi$, as is easily deduced by replacing in its defining sum $x$ by $\overline{x}$ and by using $\chi(\overline{x}) = \overline{\chi}(x)$ and $\mathrm{tr}(\overline{x}) = \mathrm{tr}(x)$. The independence of the choice of $\chi$ implies that $\lambda(r)$, for any $r$, is contained in the field of $l$-th roots of unities (actually, $\lambda(r)$ is even real as follows from the easily proved facts that $\rho(S)$ is unitary, symmetric and satisfies $\rho(S)^2 = 1$.) Thus $\rho$ satisfies the properties listen in the Theorem characterizing $\rho_l$ in § 2, and hence is equivalent to $\rho_l$. Indeed, by permuting the components of the vector valued function $\xi_c$ occurring in the definition of $\rho_l$ and by multiplying by a suitable diagonal matrix we can even assume that $\rho = \rho_l$.

We now set $\Phi = \Phi_\chi + \Phi_{\overline{\chi}}$. The independence of the matrices $\rho(A)$ ($A \in \Gamma$) of the choice of $\chi$ then implies that the subspace spanned by the components of $\Phi$ is invariant under $\Gamma$, and that $\Phi|\omega(A) = \rho_l(A)\Phi$ for all $A \in \Gamma$. It is easily verified that for all $x \in \mathbb{F}(l^2)$ one has

$$\Phi(x) \in \{0, -1, 2\}^{l-1}, \quad (r\text{-th entry of } \Phi(x)) \bmod l = \begin{cases} \operatorname{tr}(x^{(l^2-1)/3}) & \text{if } \operatorname{n}(x) = r, \\ 0 & \text{otherwise,} \end{cases}$$

where $r$ runs from 1 to $l-1$.

Next we describe lifts of $(\mathbb{F}(l^2), \operatorname{n}(x)/l)$. Let $V$ be the quaternion algebra over $\mathbb{Q}$ ramified at $l$ and $\infty$. If we set $K = \mathbb{Q}(\sqrt{-l})$ then $V$ can be described as $V = K + Ku$, where $u^2 = -1/3$ and $\alpha u = u\overline{\alpha}$ for all $\alpha \in K$. The map $c = \alpha + \beta u \mapsto \overline{c} := \alpha - u\beta$ defines an anti-involution of $V$. The reduced norm $\operatorname{n}(c)$ and reduced trace $\operatorname{tr}(c)$ of a $c \in V$ are given by

$$\operatorname{n}(c) = c\overline{c} = |\alpha|^2 + \frac{1}{3}|\beta|^2, \qquad \operatorname{tr}(c) = c + \overline{c} = \alpha + \overline{\alpha}.$$

Let $\mathfrak{o}$ be the ring of integers in $K$. Note that the rational prime 3 splits in $K$ since $l \equiv -1 \bmod 3$. i.e. $3 = \mathfrak{p}\overline{\mathfrak{p}}$ with a prime ideal $\mathfrak{p}$ in $K$. (Indeed, one can take $\mathfrak{p} = 3\mathfrak{o} + (1 + \sqrt{-l})\mathfrak{o}$.) We set

$$\mathfrak{O} = \mathfrak{o} + \mathfrak{p}v, \qquad v = \begin{cases} u & \text{for } l \equiv 3 \bmod 4 \\ \frac{1+u}{2} & \text{for } l \equiv 1 \bmod 4 \end{cases}.$$

It can be easily checked that $\mathfrak{O}$ is an order in $V$ (i.e. a subring which, viewed as $\mathbb{Z}$-module, is free of rank 4). In fact, $\mathfrak{O}$ is even a maximal order since the determinant of the Gram matrix $(\operatorname{tr}(e_j\overline{e}_k))$, for any $\mathbb{Z}$-basis $e_j$ of $\mathfrak{O}$, equals $l^2$ (cf. [V, Chap. III, Corollaire 5.3].).

We now have

**Lemma.** *(1) The dual lattice of $\sqrt{-l}\,\mathfrak{O}$ w.r.t. the quadratic form $\operatorname{n}(c)/l$ is $\mathfrak{O}$. The quotient ring $\mathfrak{O}/\sqrt{-l}\mathfrak{O}$ is the field with $l^2$ elements, and the anti-involution $c \mapsto \overline{c}$ on $\mathfrak{O}$ induces the Frobenius automorphism $x \mapsto x^l$ on $\mathfrak{O}/\sqrt{-l}\mathfrak{O}$.*

*(2) Let $I \subset \mathfrak{O}$ be an $\mathfrak{O}$-left ideal, and let $n = \operatorname{n}(I)$ be the reduced norm of $I$ (i.e. the g.c.d. of the integers $\operatorname{n}(x)$ where $x$ runs through $I$). Then the dual lattice of $\sqrt{-l}I$ with respect to $\operatorname{n}(c)/ln$ is $I$. There exists a $c_0 \in I$ such that $n(c_0)/n \equiv 1 \bmod l$, and for any such $c_0$ the map $c \mapsto c_0 \cdot c$ defines an isomorphism of quadratic modules $(\mathfrak{O}/\sqrt{-l}\mathfrak{O}, \operatorname{n}(x)/l) \xrightarrow{\approx} (I/\sqrt{-l}I, \operatorname{n}(x)/ln)$.*

Here, for convenience, we use the same symbols $\overline{n}(x)/nl$ for the quadratic form on $I$ as well as for the quadratic form induced by it on $I/\sqrt{-l}I$. The Lemma follows easily from standard facts in the theory of quaternion algebras; for the reader's convenience we sketch the proof in the Appendix to this section.

The Lemma provides us with lifts $(I, n(x)/nl)$ of $(\mathbb{F}(l^2), n(x)/l)$, and we now can write down explicitly elements of $M_k(\rho_l)$.

To this end let $\Phi \colon \mathbb{F}(l^2) = \mathfrak{O}/\sqrt{-l}\mathfrak{O} \to \{-1, 0, 2\}^{l-1}$ be defined as above. Let $I$ be an $\mathfrak{O}$-left ideal, choose $c_0$ as in the Lemma, and let

$$\pi : I \to \mathfrak{O}/\sqrt{-l}\mathfrak{O}, \qquad \pi(c) = \lambda + \sqrt{-l}\mathfrak{O} \quad \text{with } c \equiv \lambda c_0 \bmod \sqrt{-l}\mathfrak{O}.$$

Finally, let $p$ be a homogeneous spherical polynomial function on $V$. If we write polynomial functions on $V$ as polynomials $p$ in $\alpha$, $\overline{\alpha}$, $\beta$, $\overline{\beta}$, then it is spherical of degree $\nu$ (with respect to any nonzero multiple of $n(c)$) if and only if $p$ is homogeneous of degree $\nu$ and satisfies

$$\left(\frac{\partial^2}{\partial\alpha\partial\overline{\alpha}} + 3\frac{\partial^2}{\partial\beta\partial\overline{\beta}}\right)p = 0.$$

Set

$$\theta(\tau; I, p) = \sum_{c \in I} \Phi(\pi(c)) \, p(c) \, q^{n(c)/n(I)l}.$$

We suppress the dependence of this function on $c_0$ since a different choice results only in multiplying $\theta(\tau; I, p)$ by a scalar. By the Theorem on Realization by theta series we then have

$$\theta(\tau; I, p) \in M_{2+\deg(p)}(\rho_l).$$

It is easy to compute these functions with aid of a computer. In fact, by a computer calculation we found

**Theorem.** *Let $l$, $k$ be as in Table 2 (in § 2). Then the space $M_k(\rho_l)$ is spanned by the series $\theta(\tau; I, p)$, where $I = \mathfrak{O}$ for $l \neq 17$, and $I = \mathfrak{O}, \mathfrak{Op}$ for $l = 17$, and where $p$ runs through the homogeneous polynomial functions on the quaternion algebra $V$ of degree $k - 2$ which are spherical with respect to the quadratic form $\mathrm{n}(c)$.*

It is an open question whether the spaces $M_k(\rho_l)$, for arbitrary $k$ or primes $l$ ($\equiv -1 \bmod l$), are always spanned by theta series of the form $\theta(\tau; I, p)$, or, more generally, which spaces $M_k(\rho)$ of vector valued modular forms at all can be generated by theta series.

As explained in section 2 we are especially interested in the one-dimensional subspace $M_k^{(\delta)}(\rho_l)$ of functions in $M_k(\rho_l)$ which are $\mathcal{O}(q^\delta)$ with $\delta$ as in Table 2. Here we have

**Theorem (Theta formulas for conformal characters).** *(1) Let $c$, $l$, $k$ and $\delta$ be as in Table 2 (in § 2), and let $I = \mathfrak{O}$ for $l \neq 17$ and $I = \mathfrak{Op}$ for $l = 17$. Then there exists a homogeneous spherical polynomial function $p$ of degree $k - 2$ such that the $\theta(\tau; I, p)$ is nonzero and satisfies $\theta(\tau; I, p) = \mathcal{O}(q^\delta)$.*

*(2) Moreover, for any $p$ with this property, there exists a nonzero constant $\kappa$ such that the components of the Fourier coefficients of $\kappa\theta(\tau; I, p)$ are rational integers. In particular, the components of $\kappa\eta(\tau)^{-2k}\theta(\tau; I, p)$ satisfy the properties of conformal characters (1) to (5) stated in § 2.*

*Proof.* (1) The existence of a $p$ with Fourier development starting at $q^\delta$ follows from the preceding theorem and the fact that the subspace $M_k^{(\delta)}(\rho_l)$ contains nonzero elements. For the latter cf. the discussion in § 2; of course, it can also be checked by a straight forward calculation using the $\theta(\tau; I, p)$ that $M_k^{(\delta)}(\rho_l)$ is one-dimensional.

(2) Because of the latter it is clear then that, for proving the second statement of the theorem, it suffices to prove that, for at least one $p$ satisfying the condition $\theta(\tau; I, p) = \mathcal{O}(q^\delta)$, the function $\theta(\tau; I, p)$ has rational Fourier coefficients.

For proving this let $P_\nu(F)$ be the set of spherical homogeneous functions $p$ on $V$ of degree $\nu$ which are defined over the subfield $F \subset \mathbb{C}$. By the latter we mean that the coefficients of $p(c)$, when written as polynomial in the coefficients of $c$ with respect to a fixed basis of $V$, are in $F$. Note that this property does not depend on the choice of the $\mathbb{Q}$-basis of $V$. Since $P_\nu(F)$ is the kernel of a differential operator which has constant rational coefficients, when written with respect to any $\mathbb{Q}$-basis of $V$, it is clear that $P_\nu(\mathbb{C}) = P_\nu(\mathbb{Q}) \otimes \mathbb{C}$, i.e. we can find a basis of $P_\nu(\mathbb{C})$ which is contained in $P_\nu(\mathbb{Q})$. But then we deduce, using the preceding theorem, that $M_k(\rho_l)$ has a basis $\theta_j$ ($1 \leq j \leq d$) whose Fourier coefficients $a_{\theta_j}(r)$ ($r = 1, 2, \ldots$) are elements of $\mathbb{Q}^{l-1}$. For deducing this note that $\theta(\tau; I, p)$ for $p \in P_\nu(\mathbb{Q})$ has rational Fourier coefficients since $\Phi(x)$ is rational. The elements of $M_k^{(\delta)}(\rho_l)$ are now the linear combinations $\sum_j c_j \theta_j$ such that $\sum_j c_j a_{\theta_j}(r) = 0$ for all $1 \leq r < l\delta$. Since the latter system of linear equations is defined over $\mathbb{Q}$ and has a nonzero solution by part (1) we conclude the existence of rational nonzero solution, i.e. the existence of a linear combination of the $\theta_j$ with Fourier coefficients in $\mathbb{Q}$. □

If we pick a $p$ as described in the theorem, and if we denote by $\xi_{c,h}$ the $r$-th component of $\eta^{-2k}\theta(\tau; I, p)$, where $\frac{r}{l} - \frac{k}{12} \equiv h - \frac{c}{24}$ mod $\mathbb{Z}$ then it is clear that these functions satisfy properties (1) to (5) of § 2 (after multiplied by a constant, if necessary). Hence, by the uniqueness result cited in section 1, they are up to a constant the conformal characters of the $\mathcal{W}$-algebras introduced in the same section. In fact, the $\xi_{c,h}$ ($l \neq 5$) have interesting product expansions, which we shall discuss elsewhere; from these product expansions it can immediately read off that they can be normalized such that their Fourier coefficients are even non-negative integers, as its should be for conformal characters.

**Appendix.**

*Proof of the Lemma.* (1) Let $l \equiv -1$ mod 4. For $c = \alpha + u\beta \in V$ we have

$$\mathrm{tr}(c\overline{\mathfrak{O}}\sqrt{-l})/l = \mathrm{tr}(\alpha\,\mathfrak{o}/\sqrt{-l}) + \mathrm{tr}(\beta\,\overline{\mathfrak{p}}/3\sqrt{-l}).$$

Thus the left hand side is in $\mathbb{Z}$ if and only if each of the two terms on the right are in $\mathbb{Z}$. The latter is easily checked to be equivalent to $\alpha \in \mathfrak{o}$ and $\beta \in \mathfrak{p}$, i.e. to $c \in \mathfrak{O}$. The case $l \equiv 1$ mod 4 can be treated similarly, and is left to the reader.

It is clear that $\mathfrak{O}/\sqrt{-l}\mathfrak{O}$ is a ring of characteristic $l$ with $l^2$ elements. Hence it is isomorphic to a ring extension of $\mathbb{F}(l) = \mathbb{Z}/l\mathbb{Z}$ with $l^2$ elements. Moreover, it contains a root of $X^2 + 3$, namely $3u + \sqrt{-l}\mathfrak{O}$. Since $-3$ is not a quadratic residue modulo $l$ the polynomial $X^2 + 3$ is irreducible over $\mathbb{F}(l)$, hence $\mathfrak{O}/\sqrt{-l}\mathfrak{O}$ is a field. The anti-involution $c \mapsto \overline{c}$ induces an automorphism of the field $\mathfrak{O}/\sqrt{-l}\mathfrak{O}$ which is nontrivial since it maps $u$ to $-u$, and which hence is the Frobenius automorphism.

(2) If $I$ is an $\mathfrak{O}$-left ideal then $I^* = \overline{I} \cdot \mathfrak{O}^*/\mathrm{n}(I)$, where, for any left ideal $I$, we use

$$I^* = \{c \in V \mid \mathrm{tr}(Ic) \in \mathbb{Z}\}$$

(We were not able to find a reference for this basic formula: it can easily be proved using adelic methods. However, we shall need it only for $I = \mathfrak{O}$ or $I = \mathfrak{Op}$ (cf. the two theorems of the preceding section), and here it can be easily verified by direct computation. We omit the details for the general case.) Thus we find $\operatorname{tr}(\sqrt{-l}I\overline{c})/\operatorname{n}(I)l \in \mathbb{Z}$ if and only if $\overline{c} \in \overline{I} \cdot \mathfrak{O}^*\sqrt{-l}$. Using $\mathfrak{O}^* = \mathfrak{O}/\sqrt{-l}$, as follows from part (1), we find that the latter statement is indeed equivalent to $c \in I$.

Left-multiplication in the quaternion algebra induces on $I/\sqrt{-l}I$ a structure of a one-dimensional $\mathfrak{O}/\sqrt{-l}\mathfrak{O}$-vector space. Let $c_0 + \sqrt{-l}I$ be a basis element. Clearly $\operatorname{n}(c_0)/\operatorname{n}(I)$ is not divisible by $l$ since otherwise $\operatorname{n}(c)/\operatorname{n}(I)$ would be divisible by $l$ for any $c \in I$ contradicting the definition of $\operatorname{n}(I)$ as g.c.d. of all $\operatorname{n}(c)$ ($c \in I$). Thus we can choose a $\lambda \in \mathfrak{O}$ with $\operatorname{n}(\lambda)\operatorname{n}(c_0)/\operatorname{n}(I) \bmod l$. Replacing $c_0$ by $\lambda c_0$ it is then clear that $c \mapsto cc_0$ induces the claimed isomorphism. □

## 5. Comparison to formulas derivable from the representation theory of Kac-Moody and Casimir $\mathcal{W}$-algebra

In this section we compare our explicit formulas for the conformal characters with the ones obtained from the representation theory of Casimir $\mathcal{W}$-algebras [F-K-W], the Virasoro algebra [R-C] and Kac-Moody algebras [K].

The last three rational models in Table 2 are minimal models of so-called Casimir-$\mathcal{W}$-algebras. For this kind of algebras the minimal models have been determined (assuming a certain conjecture) in [F-K-W]. The representation theory of the two composite rational models ($\mathcal{W}_{\mathcal{G}_2}(2, 1^{14})$ and $\mathcal{W}_{\mathcal{F}_4}(2, 1^{26})$) is well-known [R-C,K].

In order to give the explicit formulas for the conformal characters of the minimal models of the Casimir-$\mathcal{W}$-algebras, the Virasoro algebra and Kac-Moody algebras we have to fix some notation first.

Let $\mathcal{K}$ be a simple complex Lie algebra of rank $l$ and dimension $n$, $h$ ($\check{h}$) its (dual) Coxeter number, $\rho$ ($\check{\rho}$) the sum of its (dual) fundamental weights, $W$ the Weyl group and $\Lambda$ ($\check{\Lambda}$) the (dual) weight lattice of $\mathcal{K}$. For $\lambda \in \Lambda$ let $\pi_\lambda$ denote the highest weight representation with highest weight $\lambda$.

Firstly, consider the case of the three rational models of the Casimir-$\mathcal{W}$-algebras. Formulas for the central charge, the conformal dimensions and the conformal characters of rational models of Casimir $\mathcal{W}$-algebras have been derived assuming a certain conjecture [F-K-W, p. 320]:

$$c = c(p,q) = l - \frac{12}{pq}(q\rho - p\check{\rho})^2,$$

$$h_{\lambda,\check{\nu}} = \frac{1}{2pq}\left((q(\rho + \lambda) - p(\check{\rho} + \check{\nu}))^2 - (q\rho - p\check{\rho})^2\right),$$

$$\chi_{\lambda,\check{\nu}}(q) = \eta(q)^{-l} \sum_{w \in W} \sum_{t \in \check{\Lambda}} \epsilon(w) q^{\frac{1}{2pq}(qw(\lambda+\rho) - p(\check{\nu}+\check{\rho}) + pqt)^2}$$

where $p, q$ are coprime integers satisfying $\check{h} \leq p, h \leq q$ and where $\lambda \in \Lambda$ and $\check{\nu} \in \check{\Lambda}$ are so-called dominant integral weights (the abstract form of the finitely many $\lambda$, $\check{\nu}$ is described in [F-K-W]; their explicit form can e.g. be found in Appendix D of [B-E-H$^3$]).

Using these formulas for $\mathcal{B}_2$ with $c(p,q) = c(11,6) = -\frac{444}{11}$ and for $\mathcal{G}_2$ with $c(p,q) = c(17,12) = -\frac{1420}{17}$ for $\mathcal{W}(2,4)$ and $\mathcal{W}(2,6)$, respectively, one obtains the conformal characters given in the last section (as can be checked by simply comparing enough Fourier coefficients). The last rational model, of type $\mathcal{W}(2,8)$, is a rational model of $\mathcal{W}E_7$ with $c(p,q) = c(18,23)$. However, in this case the above formula for the corresponding conformal characters contains a sum over a rank 7 lattice (the dual weight lattice) and a sum over the Weyl group of $\mathcal{E}_7$ which has order 2.903.040. Therefore, this formula is of no practical use for explicit calculations in this case. However, our formula in the foregoing section involves only a sum over a rank 4 lattice which is easy to implement on a computer.

Secondly, consider the rational models $\mathcal{W}_{\mathcal{G}_2}(2,1^{14})$ and $\mathcal{W}_{\mathcal{F}_4}(2,1^{26})$. These rational models are "tensor products" of the Virasoro minimal model with $c = -\frac{22}{5}$ and the rational model associated to the level 1 Kac-Moody algebra of $\mathcal{G}_2$ or $\mathcal{F}_4$, respectively. The two conformal characters of the Virasoro minimal model with central charge $c = -\frac{22}{5}$ are given by [R-C]

$$\chi_0^{Vir}(q) = q^{\frac{11}{60}} \prod_{n \equiv \pm 2 \bmod 5} (1-q^n)^{-1}, \qquad \chi_{-1/5}^{Vir}(q) = q^{-\frac{1}{60}} \prod_{n \equiv \pm 1 \bmod 5} (1-q^n)^{-1}.$$

The characters of rational models associated to the level 1 Kac-Moody algebras are well known from the Kac-Weyl formula [K, p. 173]. The rational model associated to the level $k$ Kac-Moody algebra of $\mathcal{K}$ has the following central charge and conformal dimensions

$$c^{\mathcal{K}}(k) = \frac{12k}{\check{h}(\check{h}+k)} \rho^2, \qquad h_\lambda^{\mathcal{K}} = \frac{(\rho+\lambda)^2 - \rho^2}{2(\check{h}+k)} \quad ((\lambda,\psi) \leq k)$$

where $\psi$ is the highest root of $\mathcal{K}$. The corresponding characters read

$$\chi^{\mathcal{K},\lambda}(q) = \eta(q)^{-n} q^{\frac{n-c^{\mathcal{K}}(k)}{24}} \sum_{t \in \check{\Lambda}} \dim(\pi_{\rho+\lambda+(\check{h}+k)t}) q^{\frac{(\rho+\lambda+(\check{h}+k)t)^2 - \rho^2}{2(\check{h}+k)}}.$$

The two conformal characters associated to the level 1 Kac-Moody algebras of $\mathcal{G}_2$ and $\mathcal{F}_4$ are given by:

$$\chi_0^{\mathcal{K}} = \chi^{\mathcal{K},0} \quad \chi_h^{\mathcal{K}} = \chi^{\mathcal{K},\lambda_1}$$

where $h = h_{\lambda_1}^{\mathcal{K}} = 2/5$ or $3/5$ and $\lambda_1$ is the fundamental weight of $\mathcal{G}_2$ or $\mathcal{F}_4$ with $\dim(\pi_\lambda) = 7$ or 26, respectively.

Using these formulas one obtains exactly the four conformal characters of the models $\mathcal{W}_{\mathcal{G}_2}(2,1^{14})$ and $\mathcal{W}_{\mathcal{F}_4}(2,1^{26})$:

$$\chi_0 = \chi_0^{Vir} \cdot \chi_0^{\mathcal{K}}, \quad \chi_{-1/5} = \chi_{-1/5}^{Vir} \cdot \chi_0^{\mathcal{K}}, \quad \chi_h = \chi_0^{Vir} \cdot \chi_h^{\mathcal{K}}, \quad \chi_{h-1/5} = \chi_{-1/5}^{Vir} \cdot \chi_h^{\mathcal{K}}$$

with $\mathcal{K} = \mathcal{G}_2, \mathcal{F}_4$ and $h = 2/5, 3/5$ respectively. The product formulas for the Virasoro characters and the formula for the conformal characters associated to the Kac-Moody algebras show that the Fourier coefficients of the two rational models are positive integers. Indeed, as one can show by comparing enough Fourier coefficients, these conformal characters are equal to the ones computed in the last section.

## Acknowledgments

W. E. would like to thank the research group of W. Nahm for many useful discussions. Parts of the present article were written by the second author during his membership at the Mathematical Sciences Research Institute at Berkeley, and he would like to thank the staff of the MSRI for providing such a warm and fruitful atmosphere. All computer calculations have been performed with the computer algebra package PARI-GP [GP].

## References

[B-E-H³] R. Blumenhagen, W. Eholzer, A. Honecker, K. Hornfeck, R. Hübel, *Coset Realization of Unifying W-Algebras*, preprint BONN-TH-94-11, DFTT-25/94, hepth-th/9406203.

[E-S] W. Eholzer, N.-P. Skoruppa, *Modular Invariance and Uniqueness of Conformal Characters*, preprint BONN-TH-94-16, MPI-94-67.

[F-K-W] E. Frenkel, V. Kac, M. Wakimoto, *Characters and Fusion Rules for W-Algebras via Quantized Drinfeld-Sokolov Reduction*, Commun. Math. Phys. **147** (1992), 295-328.

[GP] C. Batut, D. Bernardi, H. Cohen, M. Olivier, *PARI-GP* (1989), Université Bordeaux 1, Bordeaux.

[H-M] J. Milnor and D. Husemoller, *Symmetric bilinear forms*, Springer, Berlin-Heidelberg-New York, 1973.

[K] V. Kac, *Infinite Dimensional Lie Algebras and Groups*, World Sientific, Singapore, 1989.

[N] A. Nobs, *Die irreduziblen Darstellungen der Gruppen $SL_2(Z_p)$, insbesondere $SL_2(Z_2)$ I*, Comment. Math. Helvetici **51** (1976), 465-489.

[N-W] A. Nobs, J. Wolfart, *Die irreduziblen Darstellungen der Gruppen $SL_2(Z_p)$, insbesondere $SL_2(Z_2)$ II*, Comment. Math. Helvetici **51** (1976), 491-526.

[R-C] A. Rocha-Caridi, *Vacuum Vector Representations of the Virasoro Algebra*, in 'Vertex Operatos in Mathematics and Physics', S. Mandelstam and I.M. Singer, MSRI Publications Nr. 3, Springer, Heidelberg, 1984.

[Sch] B. Schoeneberg, *Elliptic Modular Functions*, Die Grundlehren der mathematischen Wissenschaften in Einzeldarstellungen, Bd 204, Springer, New Yory -Heidelberg - Berlin, 1974.

[V] M-F. Vign eras, *Arithmétique des algèbres de quaternions (Lecture Notes in Mathematics 800)*, Springer, Berlin, 1980.

Max-Planck-Institut für Mathematik Bonn, Gottfried-Claren-Strasse 26, 53225 Bonn, Germany
*E-mail address*: eholzer@mpim-bonn.mpg.de

Physikalisches Institut der Universität Bonn, Nussallee 12, 53115 Bonn, Germany
*E-mail address*: eholzer@avzw01.physik.uni-bonn.de

Université Bordeaux I, U.F.R. de Mathématiques et Informatique, 351 rue de la Libération, 33405 Talence, France
*E-mail address*: skoruppa@ceremab.u-bordeaux.fr